\documentclass[nofootinbib,aip,reprint,onecolumn]{revtex4-1}
\usepackage{mathptmx}
\usepackage{verbatim}
\usepackage{fdsymbol}
\usepackage{graphicx}
\usepackage{hyperref}

\begin{document}

\title{Bifurcation structure of interval maps with orbits homoclinic to a saddle-focus}

\author{Carter Hinsley}
\email{chinsley1@student.gsu.edu}
\affiliation{Department of Mathematics \& Statistics,  Georgia State University, 
25 Park Pl. NE, Atlanta, GA 30303, USA.}
\author{James Scully}
\email{jscully2@student.gsu.edu}
\affiliation{Neuroscience Institute, Georgia State University, 
100 Piedmont Ave., Atlanta, GA 30303, USA.}
\author{Andrey L. Shilnikov}
\email{ashilnikov@gsu.edu}
\affiliation{Neuroscience Institute, and Department of Mathematics \& Statistics, Georgia State University, \\
100 Piedmont Ave., Atlanta, GA 30303, USA.}

\date{\today}

\begin{abstract}
We study homoclinic bifurcations in an interval map associated with a saddle-focus of (2, 1)-type in $\mathbb{Z}_2$-symmetric systems.
Our study of this map reveals the homoclinic structure of the saddle-focus, with a bifurcation unfolding guided by the codimension-two Belyakov bifurcation.
We consider three parameters of the map, corresponding to the saddle quantity, splitting parameter, and focal frequency of the smooth saddle-focus in a neighborhood of homoclinic bifurcations.
We symbolically encode dynamics of the map in order to find stability windows and locate homoclinic bifurcation sets in a computationally efficient manner.
The organization and possible shapes of homoclinic bifurcation curves in the parameter space are examined, taking into account the symmetry and discontinuity of the map.
Sufficient conditions for stability and local symbolic constancy of the map are presented.
This study furnishes insights into the structure of homoclinic bifurcations of the saddle-focus map, furthering comprehension of low-dimensional chaotic systems.

\end{abstract}

\maketitle

\begin{figure}[b!]
  \begin{center}
    \includegraphics[width=0.4\linewidth]{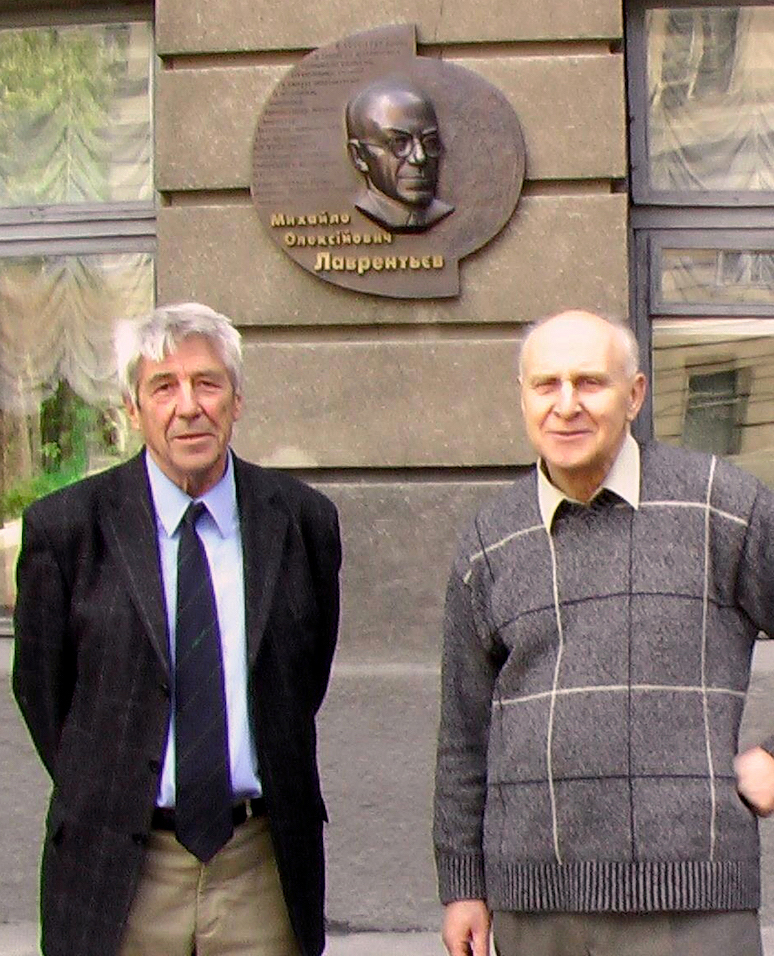}
    \caption{Two classics of high-dimensional and one-dimensional dynamics: L. P. Shilnikov and O. M. Sharkovsky (Kiev, 2005). During this visit L. P. Shilnikov was awarded the Lavrentiev medal by the National Academy of Sciences of Ukraine for his pioneering contributions to dynamical system theory.}\label{fig:fig1}
  \end{center}
\end{figure}

We begin with the acknowledgement that we are very grateful to the special editors who invited us to submit our recent research to this special issue. 
It is an honor for us to contribute to this volume of the Ukrainian Mathematical Journal, dedicated to the memory and the academic legacy of Olexander Sharkovsky, beginning with his seminal publications\cite{shark1,shark2} from the early 60s, his reference book \cite{shark4} co-authored with his students in the mid-80s of the previous century, and concluded with collection\cite{shark3} yesteryear, 2022.
One of our own, A.L.S., had the privilege of knowing Dr. Sharkovsky personally through various academic rendezvous, their initial encounter taking place at a meeting in Jurmala, Latvia in 1989. 
Predominantly, these encounters were facilitated by Yuri and Volodimir Maistrenko at their scholarly gatherings held in the serene setting of peaceful Crimea. 
Moreover, A.L.S. had a couple of occasions to interact with Dr. Sharkovsky at his parents' abode, the residence of Leonid and Ludmila Shilnikov. 
It is notable to mention that Olexander and Leonid shared an enduring friendship and academic kinship, extending over half a century, marked by mutual respect and admiration.
Each held the other's original scientific school, founded in Kiev and Nizhny Novgorod (formerly known as Gorky) respectively, in the highest esteem.

\section{Introduction}

\begin{figure}[ht!]
\begin{center}
\includegraphics[width=.4\linewidth]{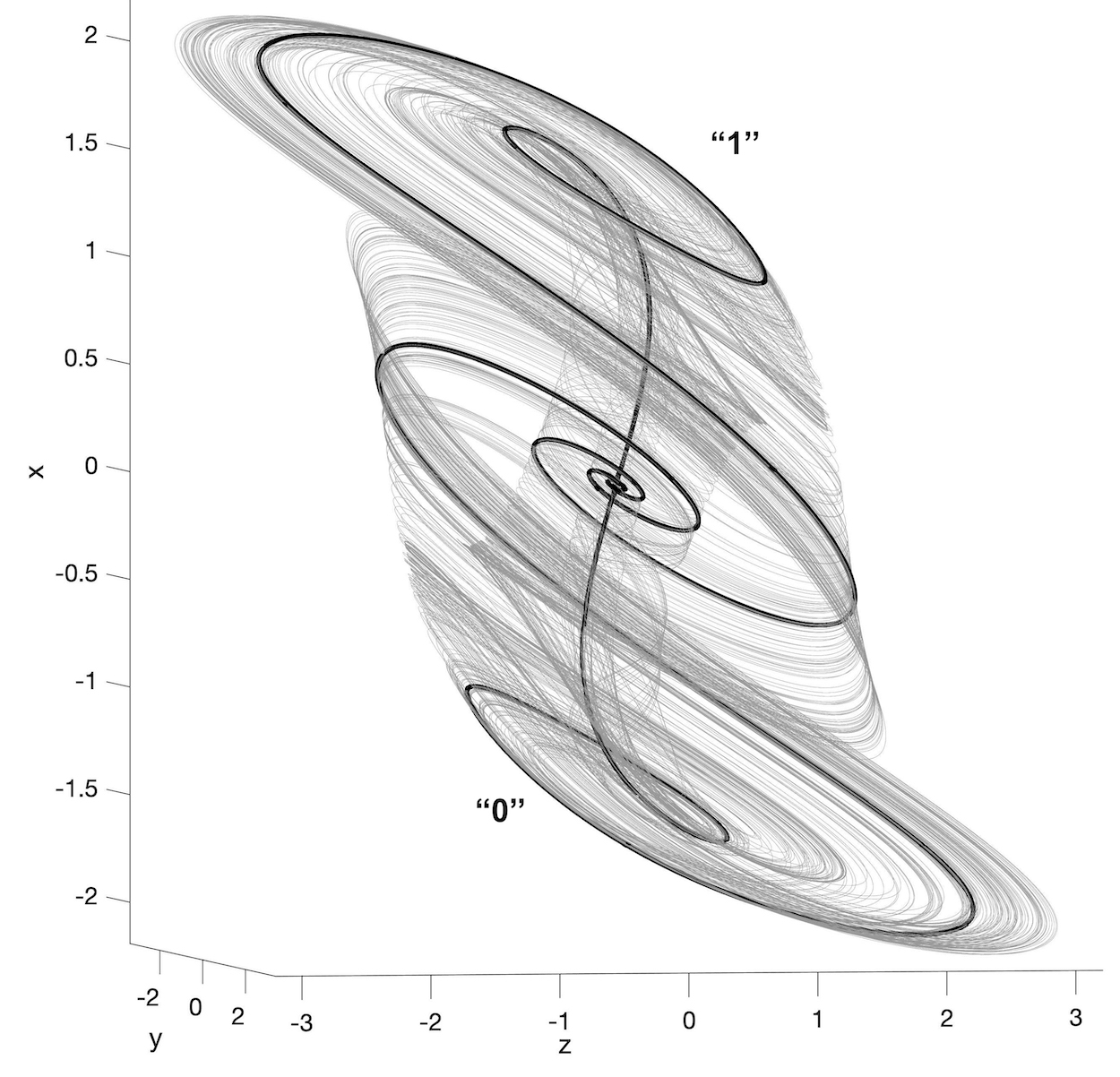}
\caption{The complex chaotic dynamics governed by the Shilnikov saddle-focus at the origin with the characteristic homoclinic figure-8 (in black) in the three-dimensional phase space of the ${\mathbf Z}_2$-symmetric model~(\ref{eqn:act1}) at $a=2.1593$, $b=0.7$, and $c=1.95$.}\label{fig:fig2}
\end{center}
\end{figure}

We aim to scrutinize and computationally illustrate the structure of bifurcation unfoldings of periodic and homoclinic orbits in one-dimensional saddle-focus return maps, especially with regards to the Shilnikov saddle-focus in the mirror-symmetric case.
These occurrences emerge near the primary figure-8 connection in a fully $\mathbb{Z}_2$-symmetric system. 
Figure~\ref{fig:fig2} offers a glimpse of such intricate dynamics, portraying the chaotic trajectories recurrently returning nearby the saddle-focus only to spiral into the three-dimensional phase space of the characteristic model\cite{ArnCoulTres1981,xing2021ordered} with reflective ${\mathbb Z}_2$-symmetry:
\begin{equation}\label{eqn:act1}
\dot{x} = y, \quad \dot{y} = z, \quad \dot{z} = -bz + cy + ax-x^3, \quad \mbox{with}\quad 
a, b, c > 0.
\end{equation}

In L. P.~Shilnikov's seminal works on the saddle focus, he convincingly demonstrated that the presence of a single homoclinic orbit of the Shilnikov saddle-focus instigates the onset of chaotic dynamics, involving a countable number of periodic orbits in the phase space of such systems.
His pioneering theories from the 1960s firmly established and underscored the critical role of homoclinic orbits within the hierarchy of deterministic chaos in its entirety \cite{ABS:1977,ABS:77b,AfrShil1983}.

Before proceeding, it seems prudent to recapitulate some fundamental elements of the Shilnikov saddle-focus theory. 
For a comprehensive understanding, one can refer to his original papers, review articles~\cite{Shilnikov1965,LP67,LP68,LP69,Shilnikov:2007,Shilnikov_heritage,gonchenko2022leonid}, and textbooks~\cite{a2001methods,arnold2013dynamical}. 
Relevant insights can also be gleaned from previous studies \cite{G83,Belyakov1,Belyakov2,Belyakov3,OSh86e,ovsyannikov1992systems,GTGN97,GLP2007,xing2021ordered,malykh2020homoclinic} that are pertinent to both the theory and the focus of this paper.

The Shilnikov saddle-focus homoclinic bifurcation serves as a fundamental and visually accessible example of chaotic dynamics within low-dimensional systems of differential equations.
Requiring a mere three dimensions for depiction, its homoclinic orbit and adjacent trajectories lend themselves to convenient visualization. 
Further, this structure's compatibility with one-dimensional return maps enhances its value as a paradigm for the evolution of mathematical and computational tools within the realm of chaotic systems.

\begin{figure}[thb!]
\begin{center}
\includegraphics[width=0.75\linewidth]{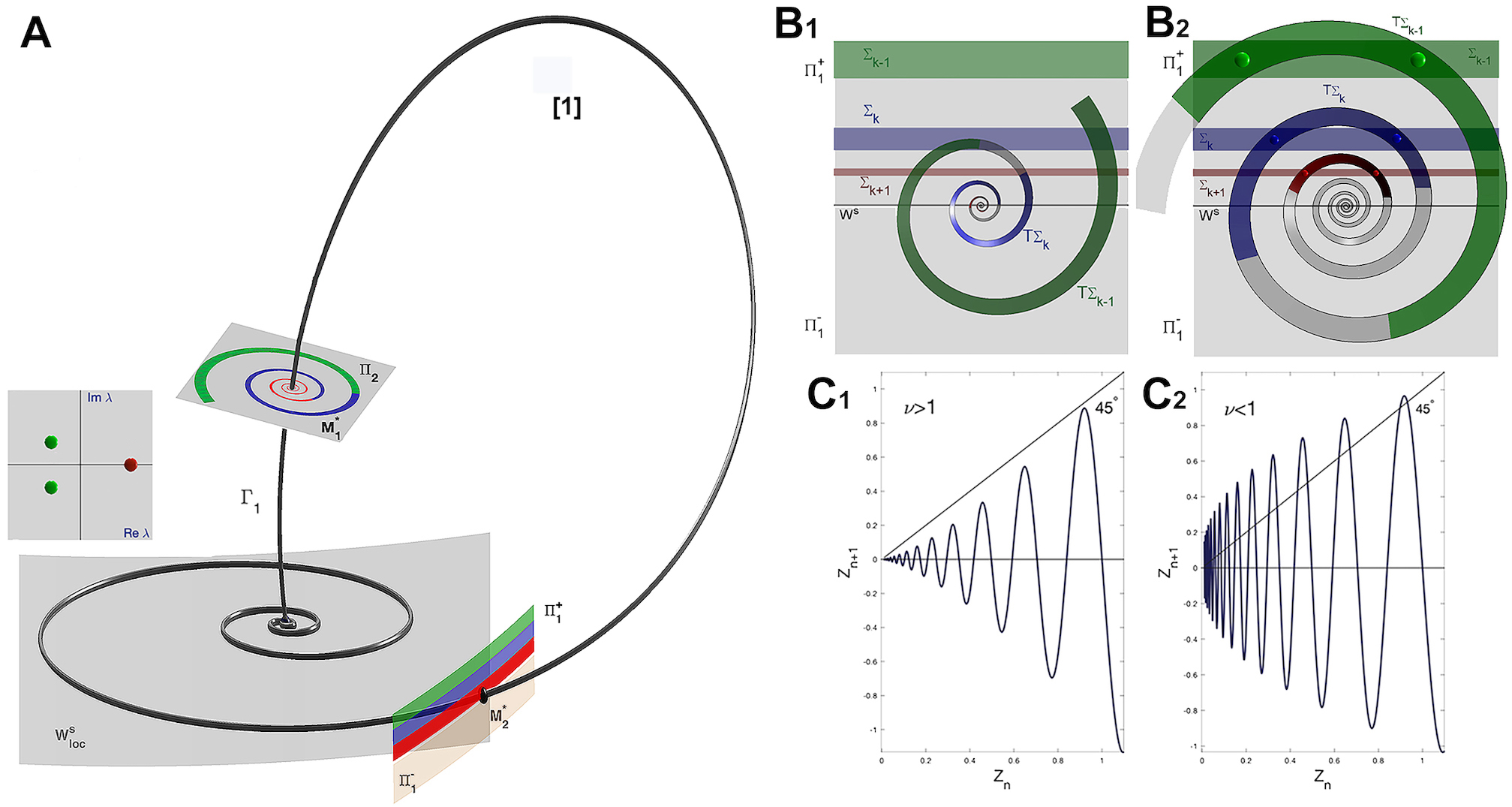}
\caption{ (A) Three-dimensional phase space showing the primary homoclinic orbit of a saddle focus of (2,1)-type, i.e., with two-dimensional stable manifold $\mathbb{W}^s$ and one-dimensional unstable manifold $\mathbb{W}^u$.  Three colored stripes painted on a two-dimensional cross-section $\Pi_1$,  locally transverse to $W^s$, are morphed along trajectories passing by the saddle focus into a colored spiral on the top cross-section $\Pi_2$, transverse to $\mathbb{W}^u$. 
 (B$_1$) The two-dimensional Poincar\'e return map $T : \Pi_1 \to \Pi_1$ is a contraction when the saddle index $\rho > 1$; the corresponding one-dimensional map is shown in C1. 
 (B$_2$)  When the Shilnikov condition  $\rho < 1$ is fulfilled, then the map is an expansion with overlapping  $T\, \Sigma_k \cap \Sigma_k$ that gives rise to countably many Smale horseshoes and saddle periodic orbits corresponding to repelling fixed points in the respective one-dimensional map in panel C$_2$; courtesy of Ref.~\cite{Shilnikov:2007}}\label{fig:shilnikovhomoclinic}
\end{center}
\end{figure}

Figure~\ref{fig:shilnikovhomoclinic}A illustrates the primary homoclinic orbit to a saddle-focus of the differential (2,1)-type. 
The designation (2,1)-type implies that the saddle-focus possesses a pair of complex conjugate characteristic exponents, denoted as $\lambda_{1,2} = -\alpha \pm i\omega$,  $\alpha, \omega>0$ (small green dots in the inset of fig.~\ref{fig:shilnikovhomoclinic}A), residing in the open left-half of the complex plane, alongside a single positive real exponent $\lambda_3$ (red dot).
It is important to stress that, for the Shilnikov saddle-focus classification, the complex pair should be the closest to the imaginary axis; this corresponds to chaos due to the existence of countably many saddle periodic orbits intersecting any small neighborhood of the saddle-focus.
On the other hand, if the Shilnikov condition is not met (i.e., if the real eigenvalue is closest to the imaginary axis), then there exists a neighborhood of the saddle-focus not intersecting any periodic orbits.~\cite{Shilnikov:2007}

System trajectories passing nearby the saddle-focus effectively map a local cross-section $\Pi_1^+$ (transverse to flow in the two-dimensional stable manifold $\mathbb{W}^s_{\rm loc}$) onto another cross-section $\Pi_2$ (transverse to the one-dimensional unstable separatrix $\Gamma_1$).
Consequently, three colored stripes delineated on $\Pi_1^+$ morph into a correspondingly colored spiral on $\Pi_2$.
The global map $\Pi_2 \to \Pi_1$ transposes the spiral back onto the original section as depicted in figs.\ref{fig:shilnikovhomoclinic}B$_1$ and \ref{fig:shilnikovhomoclinic}B$_2$.
The saddle index $\rho = \lambda_3/\alpha$ being less or greater than $1$ engenders two distinct outcomes of such a homoclinic bifurcation.
When $\rho > 1$ , i.e., local stability ``dominates'' local instability at the saddle-focus, the resulting two-dimensional map is a contraction (fig.~\ref{fig:shilnikovhomoclinic}B$_1$).
Its one-dimensional projection is visually represented in the Lammerey cobweb diagram presented in  fig.~\ref{fig:shilnikovhomoclinic}C$_1$, capturing the essential details of the map.
In accordance with Ref.~\cite{a2001methods}, we can adopt the following truncated form of the generic one-dimensional saddle-focus map:
\begin{equation}\label{eqn:1dmap}
x_{n+1} = \mu + x_n^\rho\cos(\omega\ln (x_n) +\phi)\quad  \mbox{with} \quad x_n \geq 0.
\end{equation}
In the $\mathbb{Z}_2$-symmetric case, the map becomes discontinuous for $\mu \neq 0$:
\begin{equation}\label{eqn:1dmaps}
x_{n+1} = {\rm sign}(x_n) \left[\mu + |x_n|^\rho\cos(\omega\ln|x_n| +\phi)\right].
\end{equation}
Note that the $x$ coordinate in this system does not correspond to $x$ in the (\ref{eqn:act1}) system.

\begin{figure*}[h]
\begin{center}
    \includegraphics[width=0.6\linewidth]{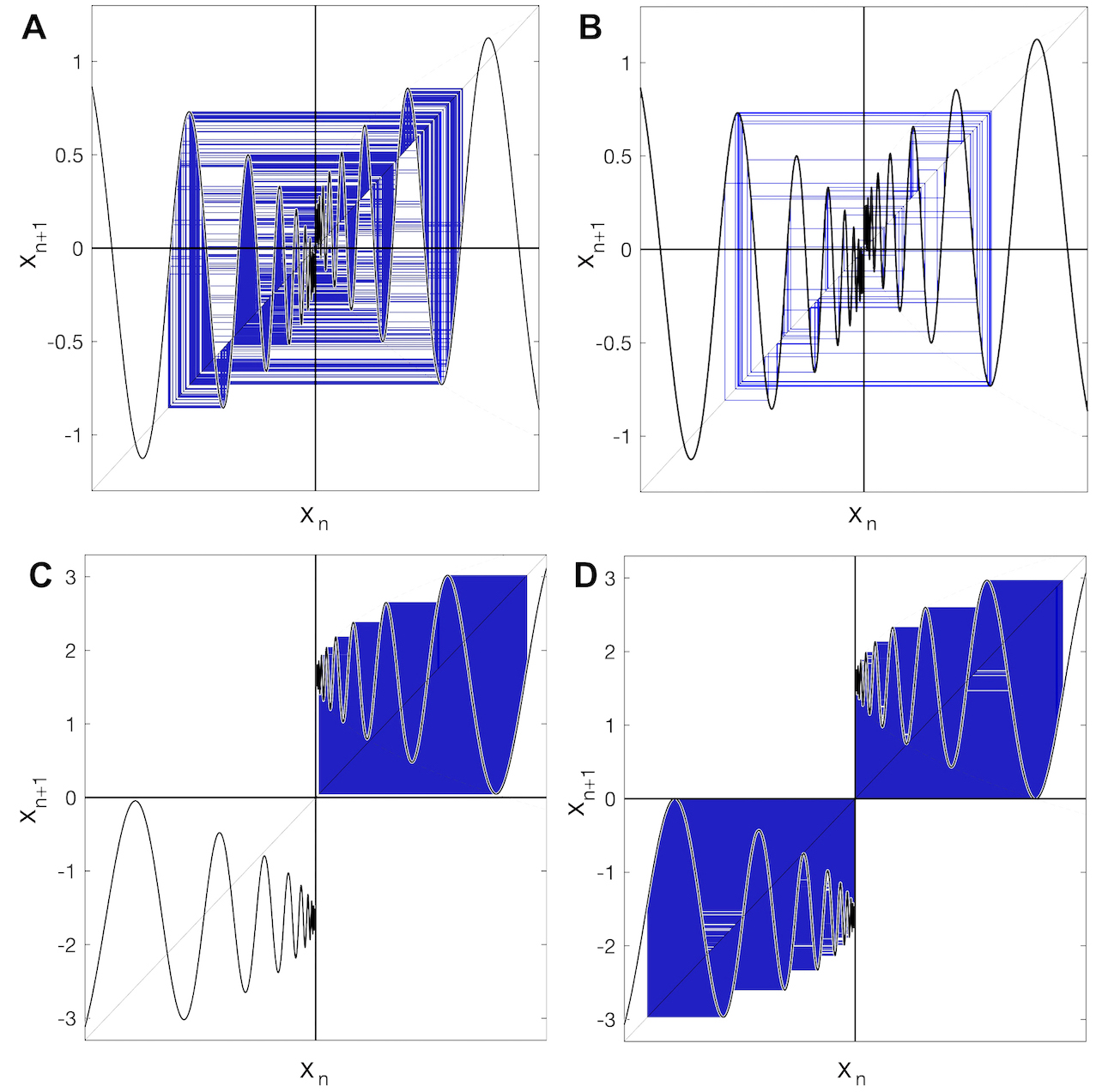}
    \caption{Snapshots of the symmetric discontinuous map~(\ref{eqn:1dmaps}) with $\rho = 0.5$ and $\omega = 10$ depicting (A) chaotic dynamics at $\mu = 0.05$, (B) a stable period-2 orbit at $\mu=0.125$, a transition from one-sided chaos at $\mu = 1.65$ in (C) to symmetric chaos at $\mu = 1.6$ in (D) after the ``boundary crisis'' when a critical point lowers below the horizontal axis.}\label{fig:4maps}
\end{center}
\end{figure*}

The parameters of this system correspond to geometric properties of the saddle-focus in the differential system: $\rho$ is the saddle index, $\omega$ is the focal frequency, and $\mu$ is the splitting parameter.
In particular, $\mu = 0$ when there is a homoclinic orbit to the saddle-focus passing once through $\Pi_1$, while $\mu \neq 0$ corresponds to the distance from the stable manifold $\mathbb{W}^s_{\rm loc}$ to the image of the origin (corresponding to the first intersection of $\Gamma_1$ with $\Pi_2$) under the map $\Pi_2 \to \Pi_1$ given by the flow.
This allows us to track the system's behavior as it undergoes a primary homoclinic bifurcation as $\mu$ crosses $0$, as well as to study secondary, tertiary, and countably many other ancillary homoclinic bifurcations of the saddle-focus as it merges with the corresponding nearby periodic orbits for $\mu \neq 0$ in the Shilnikov case $\rho < 1$.

The origin $x = 0$ in the one-dimensional map always corresponds to the saddle-focus of the three-dimensional system. For $\mu = 0$ and $\rho > 1$ (when the two-dimensional return map $T : \Pi_1 \to \Pi_1$ sends small neighborhoods of the origin into themselves) the fixed point $x^* = 0$ of the one-dimensional map (\ref{eqn:1dmaps}) is superstable.
In contrast, the scenario when $\rho < 1$ is an expansion, as depicted in fig. \ref{fig:shilnikovhomoclinic}B$_2$.
In this case, the colored (green, blue, and red) stripes do not bound or exceed their images in the expanding spiral in distance from the origin, but instead intersect their image sets.
Such intersections are interpreted as the mechanism instigating the formation of countably many Smale horseshoes, resulting in countably many unstable periodic orbits and the onset of complex dynamics in close proximity to the primary homoclinic orbit in the phase space of the differential system.
The corresponding one-dimensional return map illustrated in fig.~\ref{fig:shilnikovhomoclinic}C$_{1,2}$ locally exhibits countably many characteristic oscillations, resulting in countably many unstable fixed points at the intersections of the graph with the identity line.
It is worth mentioning that (i) these correspond to periodic orbits near the saddle-focus in the phase space of the corresponding differential system, and (ii) certain ``oscillations'' of the map graph will become tangent to the identity line as the parameters are varied, leading to new crossings or their elimination.
Such a tangency triggers a saddle-node bifurcation through which a pair of periodic orbits -- one stable and one saddle -- emerge.
It can be readily inferred that the stable orbit will soon undergo a period-doubling bifurcation when its slope in the map exceeds $1$ in absolute value; this will be succeeded by a period-doubling cascade, and so on.
This pattern is a primary reason why the Shilnikov bifurcation in three-dimensional systems is associated with the motion of the quasi-chaotic attractor \cite{AfrShil1983}, where a hyperbolic subset can coexist with stable periodic orbits emerging through saddle-node bifurcations \cite{GTS96,GST97} in a variety of models and applications \cite{BarrioShilnikov2011,BBSXS13,malykh2020homoclinic, scully2021measuring}.
This phenomenon is not necessarily observable in higher dimensions, where such homoclinic tangencies may instigate saddle-saddle bifurcations instead, as detailed in \cite{TLP98,TuraevShilnikov2008,Shilnikov_heritage}, no longer giving rise to stable periodic orbits within a chaotic attractor in the phase space.
    
In what follows we will examine the global organization of bifurcation unfoldings with biparameteric sweeps of the above one-dimensional return maps (\ref{eqn:1dmap}) and (\ref{eqn:1dmaps}) to reveal the organization of stability windows, also known as shrimps~\cite{BonattoGallas08,Gallas2010,SBU10,VitoloGlendGallas2011}, uniformly emerging in diverse applications, including models with the Shilnikov saddle focus~\cite{BarrioShilnikov2011, malykh2020homoclinic}. 

\begin{figure}[h]
\begin{center}
\includegraphics[width=.35\linewidth]{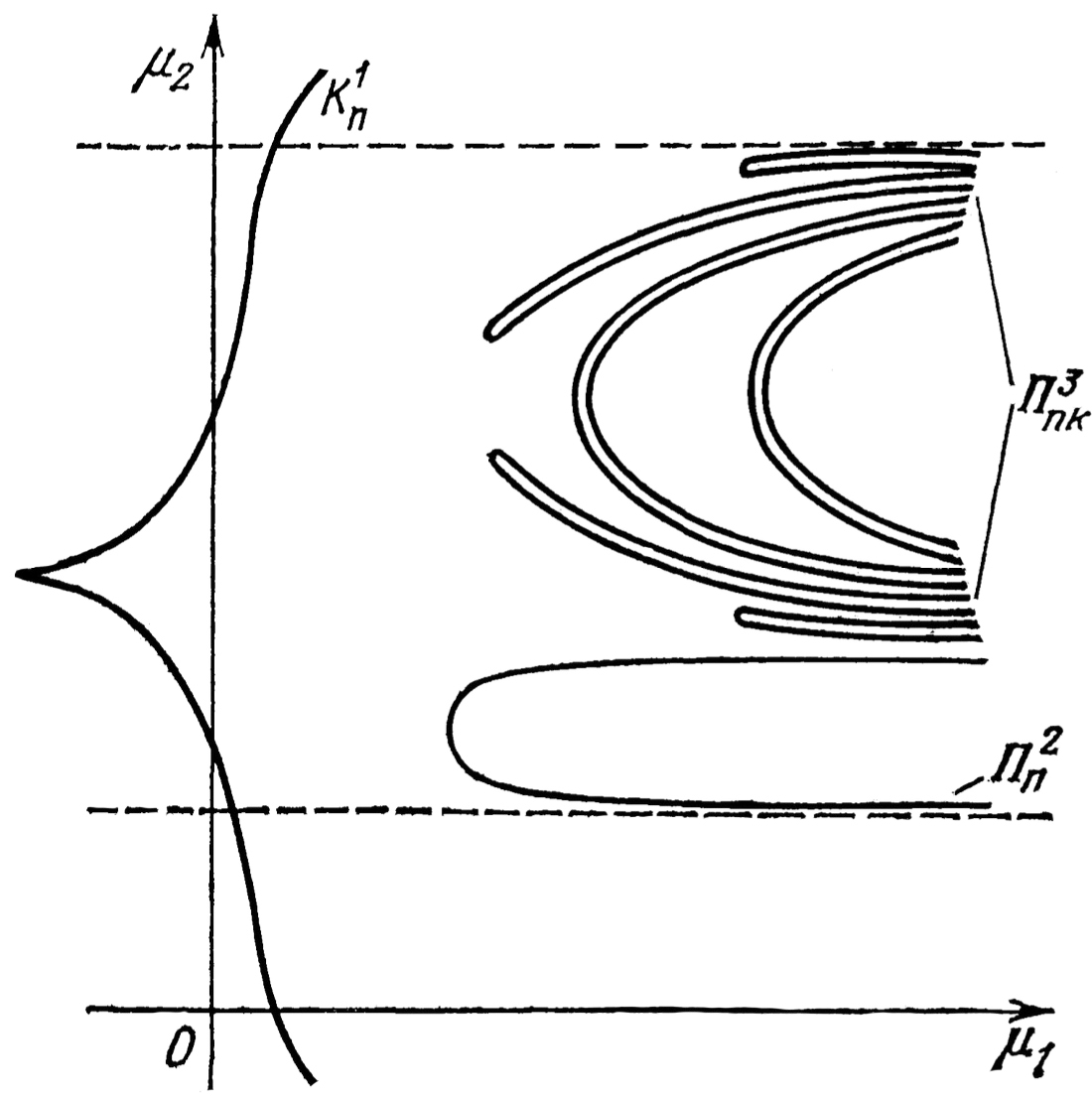}
\caption{A fragment of the Belyakov homoclinic bifurcation set\cite{Belyakov3} near the borderline transition from the Shilnikov saddle-focus for $\mu_1 > 0$ (i.e., $\rho < 1$) to a stable contraction for $\mu_1 < 0$ ($\rho > 1$).}\label{fig:belyakov}
\end{center}
\end{figure}

We will also study the fine organization of secondary and higher-order homoclinic bifurcations in such maps.
Of special consideration is the borderline codimension-2 case when the dilation map with $\rho < 1$ becomes a contraction map with $\rho > 1$.
This transition was first analytically studied by L.A.~Belyakov~\cite{Belyakov3}; see his bifurcation diagram presented in fig.~\ref{fig:belyakov}, where $\mu_1=1-\rho$, while  $\mu_2$ can be either the frequency $\omega$ or the splitting parameter $\mu$ shifting the maps given by Eqs.~(\ref{eqn:1dmap}) and (\ref{eqn:1dmaps}) up and down.
Here, a ``$\{$''-shaped curve with a cusp corresponds to two closest saddle-node or tangent bifurcations in the one-dimensional maps shown in figs.~\ref{fig:shilnikovhomoclinic}C$_1$ and C$_2$.
To the right from it, there are loci of U-shaped curves in the bifurcation diagram which correspond to secondary, tertiary, and higher-order homoclinic bifurcations in the differential system.

\begin{figure}[b]
  \centering
    \includegraphics[width=0.7\linewidth]{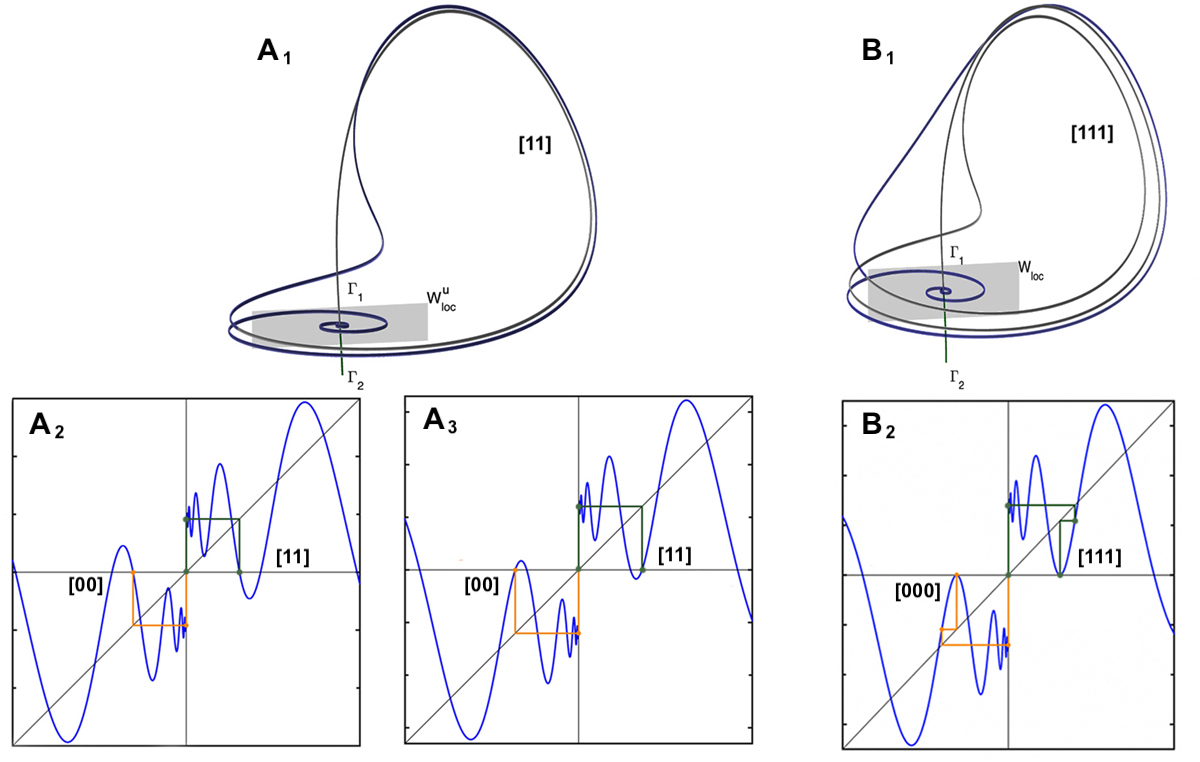}
    \caption{Comparison of {\em one-sided} secondary and tertiary homoclinic orbits. (A1) A secondary (one-sided) homoclinic orbit (coded as [11]) to the Shilnikov saddle-focus ($\rho < 1$) with its two variants in the one-dimensional return maps (A2 and A3) ending at different zeros. (B$_1$) A tertiary (one-sided) homoclinic orbit coded as [111] and its representation in the one-dimensional return map (B2) where the right forward iterates of the origin attain a critical point touching the horizontal axis -- the so-called homoclinic tangency.}\label{fig:loop11_111}
\end{figure}

To detect and differentiate such longer orbits, we employ a symbolic description, following our previous work~\cite{BSS12,xwbs, XWZS14,Xing:2014,xing2021ordered,pusuluri2018homoclinic, pusuluri2017unraveling, pusulurihomoclinicCNSNS}.
The codes [11] and [111] for the double and triple loops signify that the unstable separatrix returns to the saddle focus to complete the orbit after two and three large swings or excursions, respectively; these orbits in the differential system are secondary and tertiary homoclinics.
The respective orbits for the one-dimensional maps are demonstrated in figure panels ~\ref{fig:loop11_111}A$_2$, A$_3$, and B$_2$.
For the double loop [11] in the map~(\ref{eqn:1dmaps}), the sequence of iterates follows the pattern: $0 \mapsto \mu \mapsto 0$; whereas for the triple loop requires one more iterate: $0 \mapsto \mu \mapsto x_2 \mapsto 0$.
The oscillatory structure of the one-dimensional map allows such homoclinic orbits to emerge at different zeros or oscillatory branches as depicted in figs.~\ref{fig:loop11_111}A$_{2,3}$, though all such double orbits share the same symbolic code [11].
In Eqs.~(\ref{eqn:1dmap}) and (\ref{eqn:1dmaps}), varying $\rho$ changes the envelope of the map from convex if $\rho > 1$ to non-convex when $\rho < 1$, while the frequency parameter $\omega$ stretches and shrinks the map graph horizontally, and the splitting parameter $\mu$ shifts the graph of the one-sided map up and down.

\begin{figure}[t!]
  \centering
    \includegraphics[width=0.9\linewidth]{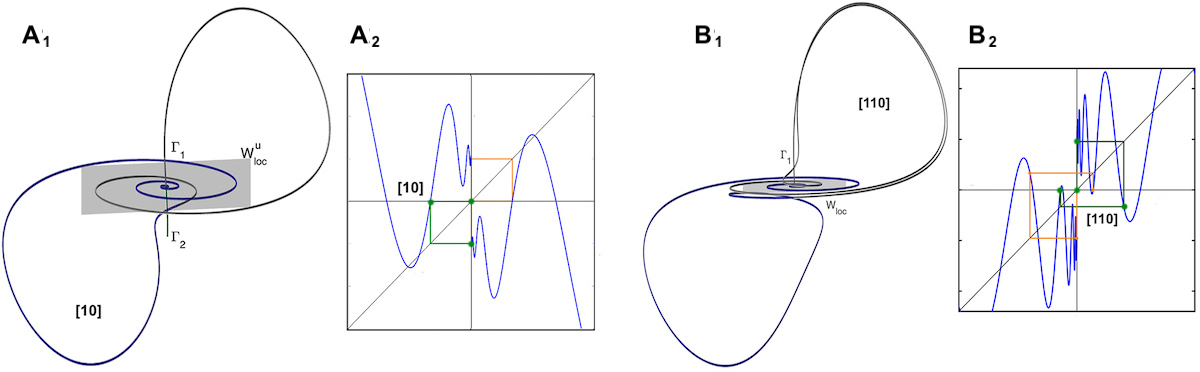}
    \caption{{\em Two-sided} secondary and tertiary homoclinic orbits and their representations in the symmetric one-dimensional map. (A1) Depicted is a secondary homoclinic orbit symbolically encoded as [10], while panel (B1) illustrates a triple homoclinic loop encoded as [110]. The corresponding return maps are displayed in (A2) and (B2) respectively.}\label{fig:loop10_110}
\end{figure}
  
We illustrate possible homoclinic orbits in a mirror-symmetric map in figs.~\ref{fig:loop10_110}A,B, showing that such orbits are inherent in $\mathbb{Z}_2$-symmetric systems like the chaotic model~(\ref{eqn:act1}) above.

\begin{figure}[b]
\begin{center}
\includegraphics[width=0.65\linewidth]{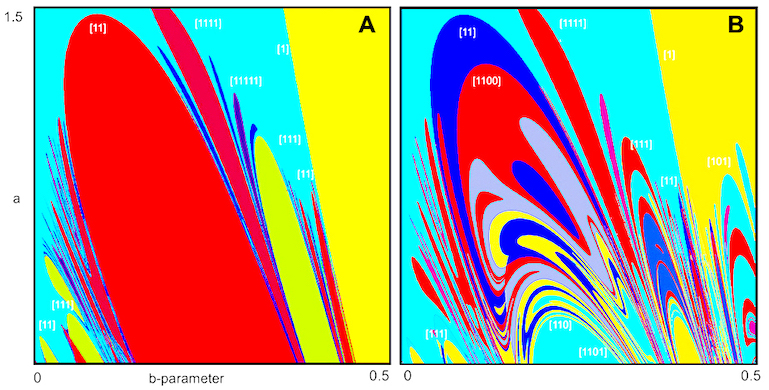}
\caption{(A) Bifurcation diagram of model~(\ref{eqn:act1}) populated by self-similar U-shaped bifurcation curves corresponding to the one-sided homoclinic orbits coded as [11], [111], ... along the boundaries of solid-color regions. (B) Fractal organization of bifurcation structures corresponding to one- and two-sided homoclinic orbits coded with all symbolic sequences.}\label{fig:u_shapes}
\end{center}
\end{figure}

This introduction concludes with snapshots showcasing the fractal organization of some global bifurcation unfolding representing a rich variety homoclinic orbits to the saddle-focus in the system~(\ref{eqn:act1}).
Figure~\ref{fig:u_shapes}A displays numerous U-shaped curves corresponding to one-sided homoclinics, while two-sided homoclinics populate within the spaces bounded by these U-shaped curves (fig.~\ref{fig:u_shapes}B). The subsequent analysis will offer a more granular examination of these structures using a computationally efficient symbolic approach.

\section{Symbolic representation and homoclinic bifurcation unfoldings}

\subsection{Partitioning the one-dimensional map}

The saddle-focus in the system corresponds to the origin $x = 0$ of the one-dimensional map, and the homoclinics in the differential system correspond to successive forward iterates of the map beginning and ending at the origin.
The map is generally discontinuous at $0$, and there are three possible behaviors at the discontinuity.
Firstly, the origin may be treated as a fixed point corresponding to the saddle focus.
The second and third possibilities involve the trajectory leaving the saddle-focus in either direction along the one-dimensional unstable manifold, corresponding to sending $0 \mapsto \mu$ and $0 \mapsto -\mu$ respectively.
We construct a binary sequence which encodes the sequence of positive and negative excursions a trajectory of the differential system takes; for each choice of parameters of the maps there correspond two such symbolic sequences.
The first element of the sequence is ``1'' for a positive excursion corresponding to $x_1 = \mu$ and ``0'' for a negative excursion corresponding to $x_1 = -\mu$.
The rest of the sequence is generated from the signs of successive iterates $x_n$, $n\ge 0$, of the chosen initial point, with ``0'' corresponding to $\rm{sign}(x_n) = -\rm{sign}(\mu)$ and ``1'' to $\rm{sign}(x_n) = \rm{sign}(\mu)$.
Due to the symmetry of the system there is a mirror image of each sequence, but we will in this paper always follow the sequence originating on the right branch of the symmetric one-dimensional map ($x_1 > 0$).

Consider the mappings from the $(\rho, \mu^+)$-parameter half-plane  to the $n^{\text{th}}$ iterates starting with the initial point $x_1 = \mu > 0$.
It is precisely the zeros of these mappings (where $x_n = 0$) that define corresponding bifurcation curves of the homoclinic orbits of the $n^{\text{th}}$ degree in the parameter space.
Reaching $x_n = 0$ is encoded symbolically as a termination of the sequence.
This sequence constitutes a binary representation of the dynamical behavior at each point, providing a comprehensive description of the homoclinic bifurcation structures.
As such, this method transforms the intricate problem of calculating homoclinic orbits in continuous-time dynamical systems into the simpler problem of finding zeros of iterates in discrete maps.
This transformation considerably simplifies the analysis and enables efficient computation of homoclinic structures.

\subsection{Basic use of the symbolic trajectory representation}

Two procedures are used to process the binary sequences.
The first procedure is to select particular sequences which illustrate particular aspects of the homoclinic structure.
The zeros of the first iterate of $\mu$ correspond to the boundary between various sequences [XX1...] and [XX0...] (here the Xs denotes various identical initial substrings in such sequences), as well as to secondary homoclinic curves in the ODE system.
Similarly, the zeros of the second iterate of $\mu$ correspond to all  bifurcation curves of tertiary homoclinic orbits.

For asymmetric systems with one-dimensional return map (\ref{eqn:1dmap}), only positive $x$ values are relevant, so the only homoclinics to consider are one-sided and correspond to sequences of repeated ``1''s.
For one-sided orbits with $\mu > 0$, it is necessary to truncate sequence just before their first zero entries.
Although in this case one cannot distinguish homoclinic orbits from non-homoclinic orbits symbolically, the boundaries of regions in parameter space corresponding to particular symbolic sequences do form homoclinic bifurcation curves.

The second procedure is to compute an embedding of binary sequences of arbitrary length into the interval $[0, 1]$.
For a binary sequence $[S_1, S_2, \ldots, S_N]$ of length N, this is computed as a partial power series with the factor $\frac12$:
\begin{equation}
K(\rho, \, \mu) = \sum_{i=1}^N S_i\frac{1}{2^i}.
\end{equation}

\subsection{Bifurcation unfoldings in the $(\rho, \mu)$-plane of the interval map}

\begin{figure*}[h]
\begin{center}
\includegraphics[width=1.0\linewidth]{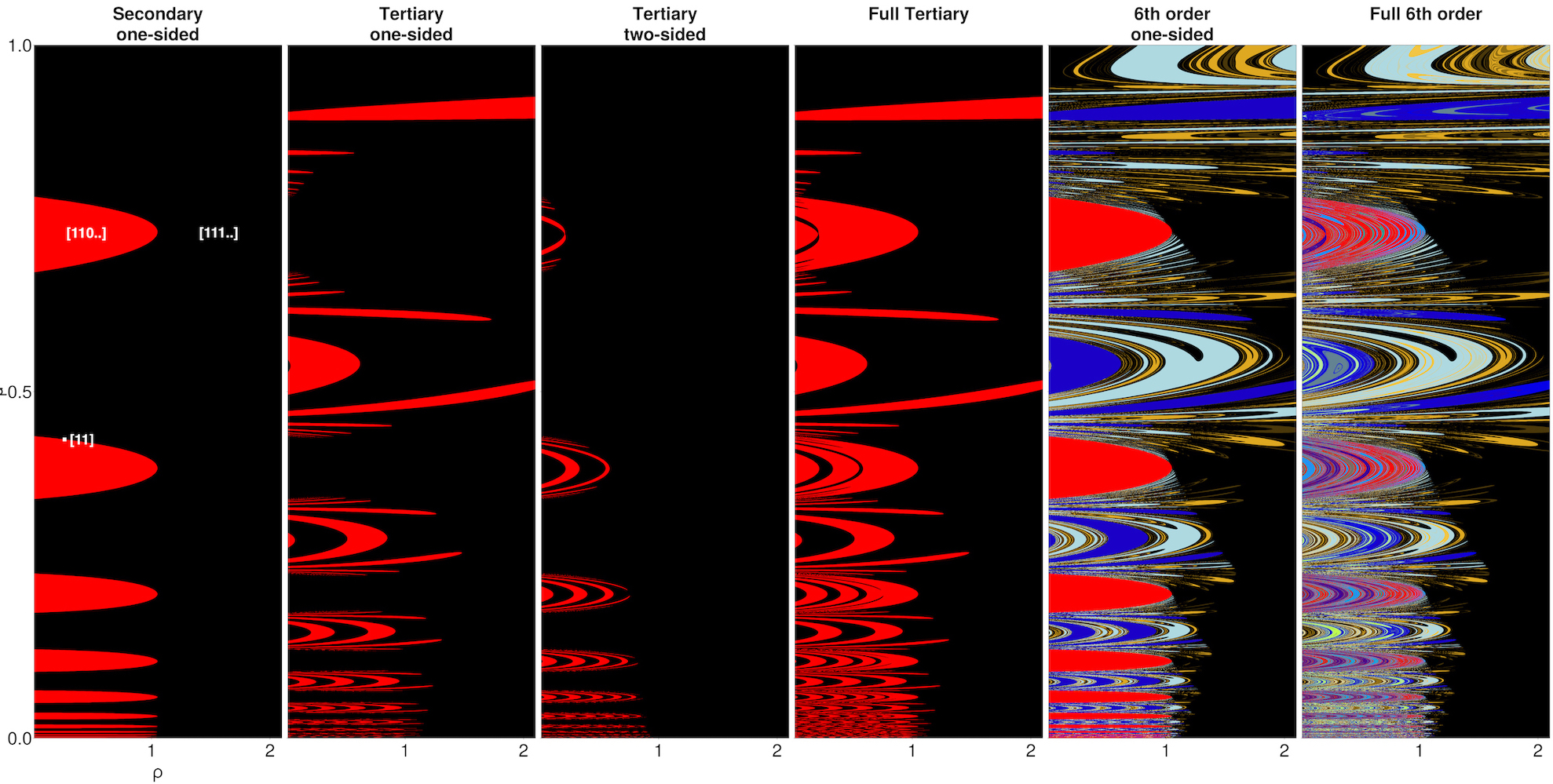}
\caption{Homoclinic bifurcation structures in the $(\rho, \mu^+)$-parameter half plane of the one-dimensional saddle-focus map. (A) The U-shaped curves of secondary homoclinics accumulating to the primary homoclinic curve at $\mu = 0$ and separating the red regions associated with sequences beginning with [110...] from the black region where sequences begin with [111...]. (B) One-sided tertiary homoclinic structure corresponding to [111]-orbits revealed by the boundary of parameter regions with sequences originating with symbols [1110...] (red), while the two-sided tertiary structure comes from separating sequences beginning with symbols [1100...] (red) and [1101...] (black) in (C). (D) The full tertiary structure corresponds to orbits differentiated by their fourth symbols, showing how the zeros of successive iterates of the initial value $\mu$ partition the parameter space along the homoclinic curves. (E, F) Bifurcation curves corresponding to the up-to-6th-order homoclinics of both types: one-sided (E) and two-sided (F) cases.}\label{fig:homoclinicstructures1}
\end{center}
\end{figure*}

The overarching structure of parameter sets for one- and two-sided sequences up to order 6 is summarized in fig.~\ref{fig:homoclinicstructures1}, with several panels presented for side-by-side comparison.
Panel~A reveals a collection of U-shaped bifurcation curves of secondary homoclinic orbits accumulating to the primary homoclinic at at $\mu = 0$ from above.
The top and bottom branches of a secondary homoclinic bifurcation curve correspond to [11]-encoded double loops occurring in the one-dimensional map as illustrated in figs.~\ref{fig:loop11_111}A$_{2,3}$: the forward iterates of the origin come back after two steps: $0 \mapsto \mu \mapsto 0$.
The peak of this U-shaped bifurcation curve at $\rho = 1$ corresponds to the case when the orbit involves a critical point of the map yields a coincidence of the graph with the horizontal axis, producing a homoclinic tangency much like the case illustrated in fig.~\ref{fig:loop11_111}B$_2$ for the tertiary homoclinic orbit.
For fixed $\rho$ and varying $\mu$ values, secondary homoclinic orbits may form at the various oscillatory branches of the one-dimensional map positioned some distances away from the origin.
This accounts for the shape and multiplicity of such U-shaped bifurcation curves, which become narrower as $\mu$ decreases, accumulating to the primary homoclinic bifurcation at $\mu=0$.
Also noteworthy is that these peaks lie exclusively on the line $\rho = 1$, with no secondary homoclinic bifurcations in the $\rho > 1$ half plane.
This implies that the secondary one-sided homoclinic tangencies are exclusive to the Shilnikov saddle-focus; i.e., where $\rho \leq 1$.
However, this is not the case for the one-sided tertiary and higher-order homoclinics, nor is it the case for two-sided homoclinic bifurcations in general, all of which will be discussed in later sections.

\begin{figure}[t]
\begin{center}
\includegraphics[width=1.\linewidth]{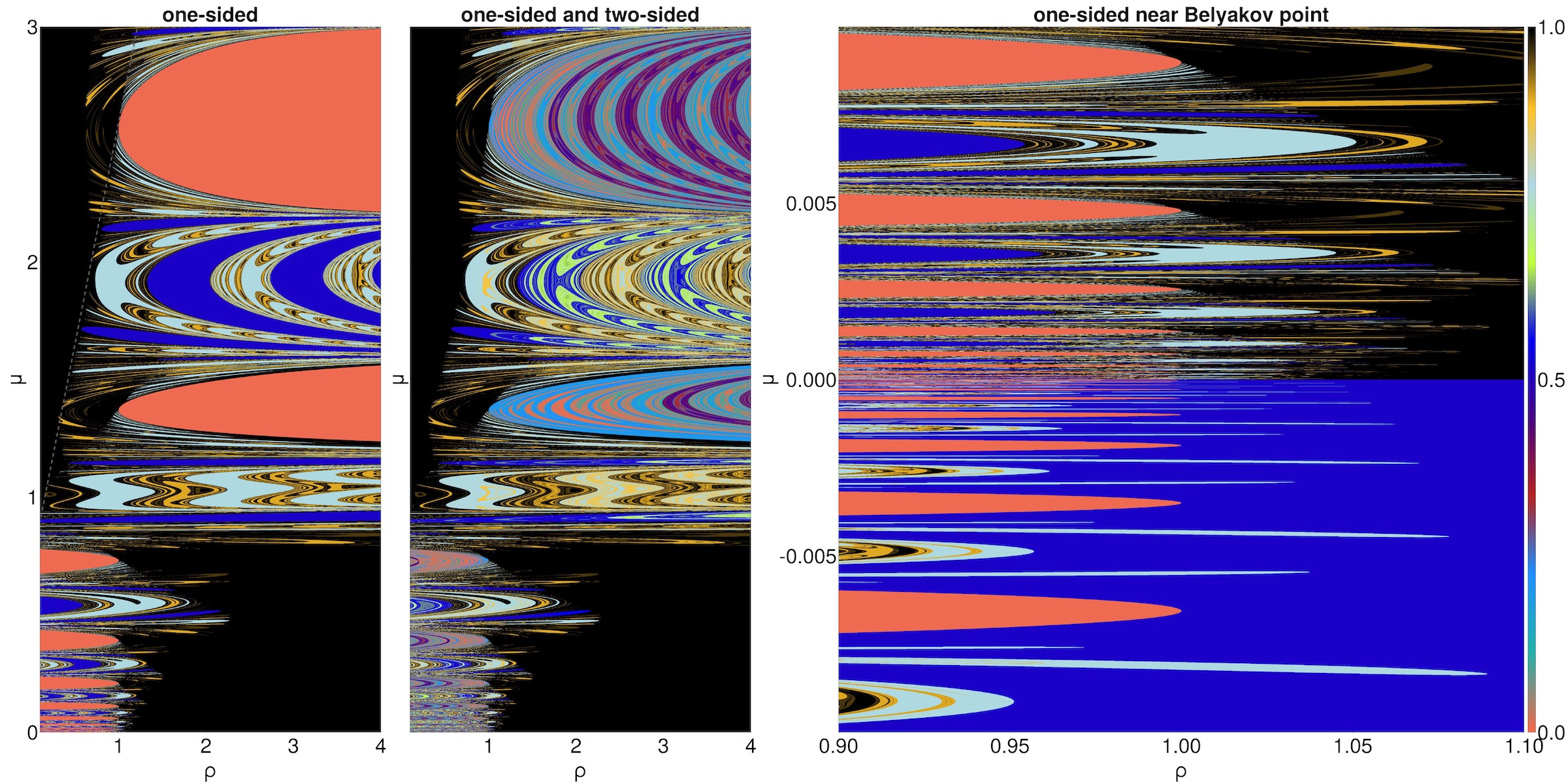}
\caption{(A, B) Homoclinic bifurcations of the one-dimensional map for large $\mu$. The left panel shows the one-sided orbits and the right panel depicts the two-sided orbits. The orientation of the homoclinic orbits in the parameter plane switches due to the changing relationship between the envelope $\mu \pm |x|^\rho$ and the image of $\mu$. (C) Low-order homoclinic bifurcations of the saddle-focus at the transitions of colors for small $\mu$ densely organized about the primary homoclinic orbit at the origin of the map at $\mu = 0$. The color bar corresponds to the embedding of a symbolic sequence at a given parameter value into the interval $[0, 1]$.
Countably many homoclinic U-shaped curves of a particular color lie tangent to each of countably many monotonic curves originating at $\rho = 1, \mu = 0$.
The U-shaped regions become uniform in this area, illustrating the fractal nature of dynamics of the saddle-focus map and its bifurcation diagram.}\label{fig:outerspace}
\end{center}
\end{figure}

While only small values of $\mu$ are relevant to the study of systems in a neighborhood of the primary homoclinic bifurcation, the behavior of the map for arbitrary $\mu$ is interesting in its own right.
In figs.~\ref{fig:outerspace}A and B, we explore the impact of larger values of $|\mu| > 1$ on homoclinic orbits.
When $|\mu|$ exceeds 1, the relationship between the envelope (due to the term $|x_n|^{\rho}$ in (\ref{eqn:1dmaps})) and the image of $\mu$ changes.
At $\mu = 1$, the envelope has a root at $x = \mu$, and thus homoclinic tangencies relevant to the flow arise only for $|x| \geq |\mu|$ so that homoclinic bifurcation curves are seen for large $\rho$ but cannot be found for $\rho$ small.
This changes the position of the homoclinic U-shaped curves, from being contained mostly within the left half of the parameter plane, to being found predominantly within the right half as depicted in these two figures.
The left panel demonstrates this effect in the case of one-sided homoclinic orbits, while the middle panel exhibits the structure of such homoclinic bifurcation curves in the two-sided case.

Figure~\ref{fig:outerspace}C demonstrates the order of homoclinic orbits and their bifurcation curves for small values of $\mu$ in the bifurcation diagram near the demarcation line $\rho = 1$ in the one-dimensional saddle-focus map, to be compared with the sketch in fig.~\ref{fig:belyakov} from the original Belyakov theory\cite{Belyakov3}.
In this case, the map exhibits fractal structure organized about the codimension-2 Belyakov point $(\rho = 1,\, \mu = 0$), with bifurcation curves of homoclinic orbits of higher orders drawn into a front at $\rho = 1$.
This observation provides an intricate look into the dynamics of the system and the fractal nature of orbits homoclinic to saddle-foci and periodic orbits in neighborhoods thereof.

\section{Stability modulation by homoclinic and shrimp structures}

\begin{figure}[t!]
\centering
\includegraphics[width=0.8\textwidth]{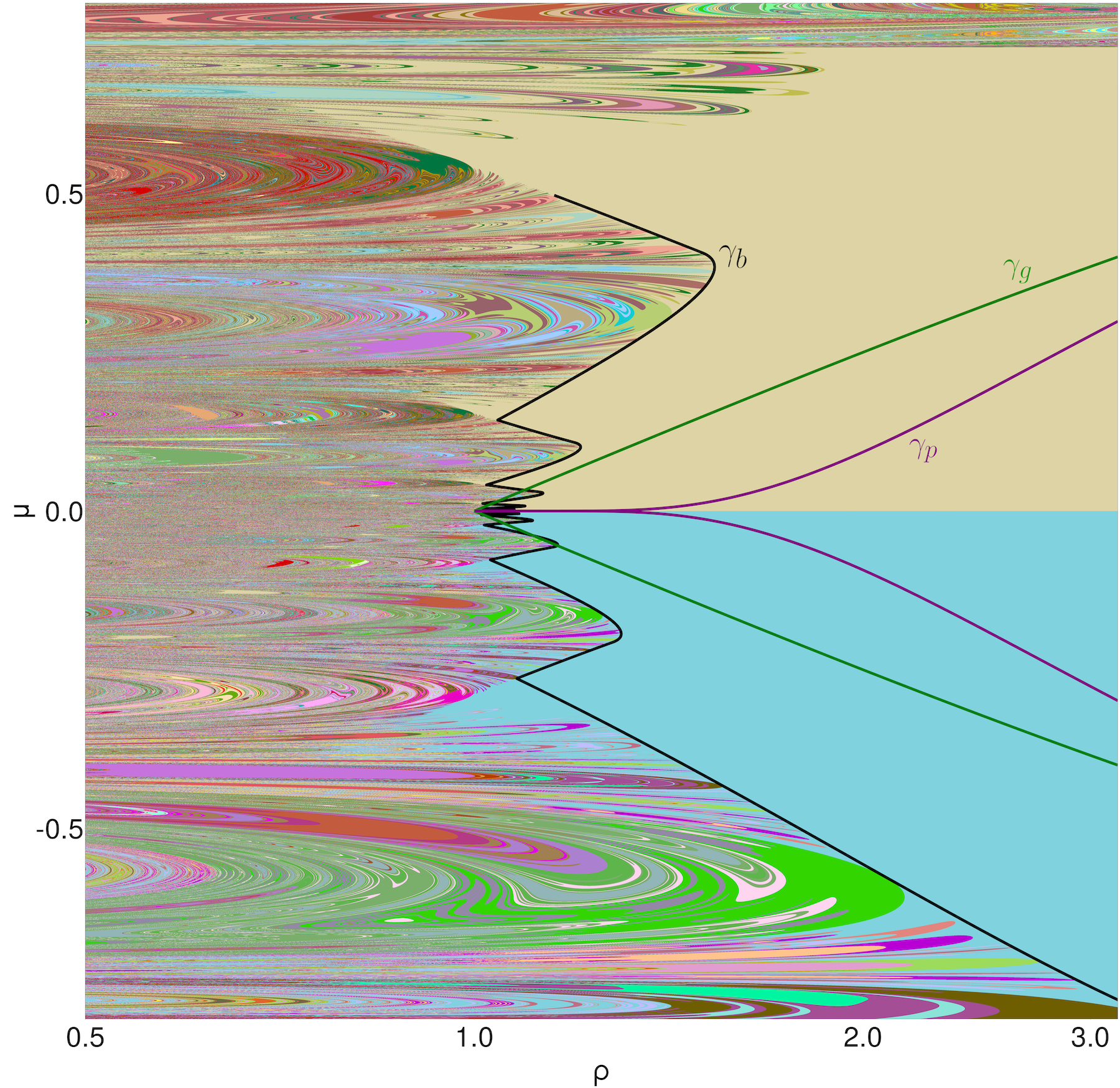}
\caption{Low- and high-order homoclinic bifurcation structure of the one-dimensional saddle-focus map with $\omega = 3.6$ on the 8000x8000 pixel scan. Outside of the cone-shaped region bounded by the green curves $\gamma_g$ for $\rho > 1$ there exists no invariant interval in the map~(\ref{eqn:1dmaps}) as its iterates may diverge for some choice of $\omega$ given $\rho, \mu$ in this region.
		The region bounded by the cusp-like purple curves $\gamma_p$ comprises parameters for which the derivative remains less than 1 in absolute value within the invariant interval, so orbits converge to the unique fixed point of the map.
		The black curves $\gamma_b$ are solutions to systems of equations corresponding to critical zeroes of the iterated map, and for sufficiently small $\mu$ serve as upper bounds on values of $\rho$ for which such ancillary homoclinic bifurcations may be found.
	}
	\label{fig:fullhomoclinicstructure}
\end{figure}

Diving deeper into the complexity of the one-dimensional discontinuous saddle-focus map~(\ref{eqn:1dmaps}), we now shift our attention to the substantial regions of parameter space known as stability windows.
It is well known that saddle-node bifurcations give rise to stability windows as various tangencies between the map graph and its higher order degrees and the identity line occur, or when its negative slope, or that of its higher degrees becomes less than one in the absolute value. These stability windows can be vividly demonstrated through the Lyapunov exponent ($LE$), which in our context is evaluated over a trajectory of 5000 iterates by taking the mean of the logarithm of the absolute derivatives of the map along the trajectory as follows:
\begin{equation}
LE (\rho, \, \mu) = \frac{1}{N} \sum_{i=1}^N \log \left | \frac{dx_{n+1}}{dx_n} (x_i) \right |.
\end{equation}

\begin{figure}[t!]
\includegraphics[width=1.0\linewidth]{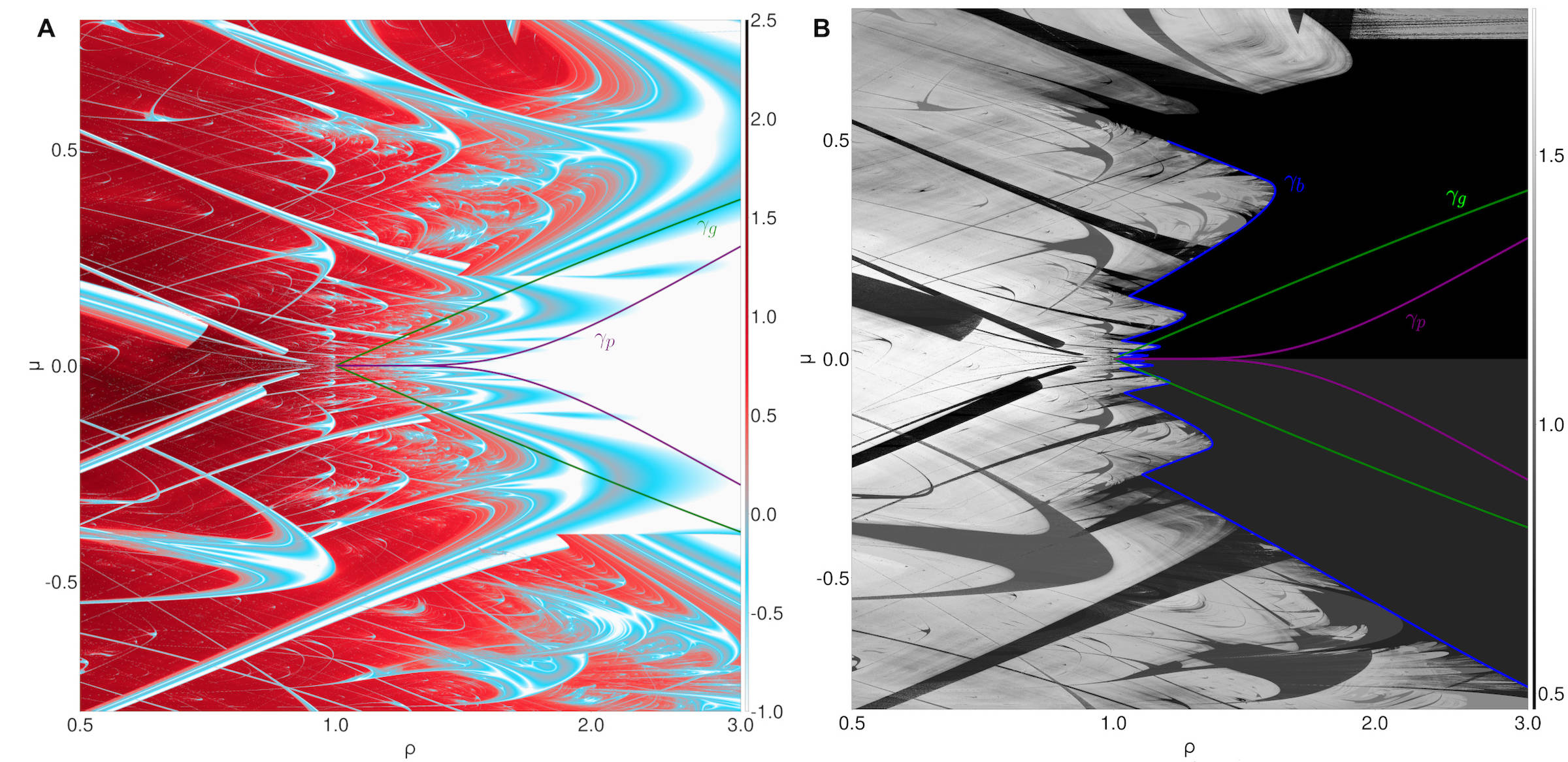}
\caption{
(A) Chaos and stability windows (``shrimps'') in a 4000x4000-pixel scan of the $(\rho, \mu)$-parameter plane of the one-dimensional saddle-focus map with $\omega=5$.
The heatmap represents the magnitudes (color bar on the right) of Lyapunov exponents computed over trajectories of length 5000 with the red color indicating chaos and the blue/white colors signifying stability.
Discontinuities in the color grading correspond to branch-switching in areas of multistability, mostly well-organized for $\rho > 1$.
(B) 4000x4000-pixel bifurcation diagram of the map viewed through the lens of Lempel-Ziv complexities of the symbolic binary sequences
of long orbits beginning with the same initial point  $x_1 = \mu$: lighter shades signify higher symbolic 
complexity (color bar on the right).
The overlaid curves $\gamma_b$, $\gamma_g$, and $\gamma_p$ are the same as in fig.~\ref{fig:fullhomoclinicstructure}. 
Despite unbounded and Lyapunov-positive trajectories in much of the $\rho > 1$ half plane, the symbolic representations of these orbits are very simple for small enough $|\mu|$, undisturbed by distant homoclinic structures.
The shrimp structures from panel~(A) appear as windows of relatively lower complexity.
}
	\label{fig:chaosmeasures}
\end{figure}

Figure~\ref{fig:chaosmeasures}A visualizes the $(\rho, \mu)$-parameter plane of the given saddle-focus map:
the color-coded heatmap reveals chaoslands in red, where $LE > 0$, and stability windows in blue and white, where $LE \le 0$.

It is worth noting that the presence of many multistability regions is a complex aspect that the Lyapunov exponent computed from a single initial value does not address directly.
The accurate and in-depth exploration and understanding of multistability principles in systems with saddle-foci remains yet an open challenge.
However, this visualization still offers insightful glimpses into the chaotic region and aids in understanding the overall stability landscape of the system.

\subsection{Stability in the absence of homoclinic interference}
It will be useful to note going forward that the derivative of the map (\ref{eqn:1dmaps}) is given by the expression
\begin{equation}\label{eqn:1dmapderivative}
\frac{\text{d}x_{n+1}}{\text{d}x_n} = \frac{|x_n|^\rho}{x_n}\sqrt{\rho^2+\omega^2}\cos\left(\omega\ln|x_n| + \tan^{-1}\frac\omega\rho\right).
\end{equation}

When the Shilnikov condition $\rho < 1$ is accompanied by the existence of a primary homoclinic (that is, when the splitting parameter $\mu = 0$ so that $x = 0$ is a fixed point of the map), chaotic behavior is observed in a neighborhood of the origin, associated with the existence of countably many unstable periodic orbits.
However, for nonzero splitting parameter $\mu$ there exist ancillary homoclinic orbits to the saddle focus, with tertiary and higher-order homoclinic orbits present even for $\rho > 1$.
The curve $\gamma_b$ in the $(\rho, \mu)$-parameter plane, seen in fig.~\ref{fig:fullhomoclinicstructure}, serves as an upper bound on $\rho$ for which homoclinic bifurcations can occur given $|\mu| \ll 1$.

$\gamma_b$ is determined in part by the explicit solution of the system of equations $\frac{\rm{d}x_{n+1}}{\rm{d}x_n}(x_1) = x_2 = 0$ for the parameter $\mu$.
This admits countably many solutions
\begin{equation}
\mu = (-1)^k\frac{\omega}{\sqrt{\rho^2+\omega^2}}\exp\left(\frac\rho\omega\left(\pi\left(-k-\frac12\right)-\tan^{-1}\frac\omega\rho\right)\right)
\end{equation}
indexed by $k \geq 0$, lying in the upper half parameter plane $\mu > 0$ for $k$ even and in the lower half plane $\mu < 0$ for $k$ odd.
The rest of $\gamma_b$ is determined by the implicit solution of the system of equations $\frac{\rm{d}x_{n+1}}{\rm{d}x_n}(x_1) = x_3 = 0$ in the $(\rho, \mu)$-plane.
Again, there are countably many solutions
\begin{align}
\begin{split}
x_2 &= (-1)^k\mu + \frac{\omega}{\sqrt{\rho^2+\omega^2}}\exp\left(\frac\rho\omega\left(\pi\left(-k-\frac12\right)-\tan^{-1}\frac\omega\rho\right)\right), \\
0 &= \mu + \left|x_2\right|^\rho\cos\left(\omega\ln\left(\left|x_2\right|\right)\right)
\end{split}
\end{align}
indexed by $k \geq 0$, this time lying in the lower half parameter plane $\mu < 0$ for $k$ even and in the upper half plane $\mu > 0$ for $k$ odd; the solution sets to these equations do belong to each half plane, but serve to bound homoclinic bifurcation sets only in one half plane or the other.
Only a certain restriction of these solution sets within the $(\rho, \mu)$-plane correspond to $\gamma_b$, although the equations involved do govern organization of homoclinic bifurcations internally to the region bounded above in $\rho$ by $\gamma_b$. Moreover, there exist conditions corresponding to higher-order iterates of the map which serve to further organize the homoclinic bifurcation structure; in general, these conditions correspond to systems of equations for which only implicit solutions may be obtained.

Through geometric analysis of the one-dimensional map, parameter values associated with the existence of a fixed point are determined. Additionally, some conditions under which bounds on trajectories can be established are identified.
As our analysis concerns behavior of the map (\ref{eqn:1dmaps}) in a small neighborhood of $x = 0$, it is useful to note that in many cases a compact invariant interval containing the origin can be given.

For $\rho > 1$ and $\mu = 0$, a small neighborhood of $x=0$ cannot contain any fixed points of the map other than the origin $x=0$ itself.
As $\frac{\text{d}x_{n+1}}{\text{d}x_n}(0) = 0$, the origin is stable.
However, for $\mu \neq 0$, orbits may wander chaotically and the non-convexity of the envelope $\left|x_{n+1} - \frac{x_n}{|x_n|}\mu\right| \leq |x_n|^\rho$ can lead to exploding trajectories.
In preventing these issues it is enough to consider only $x > 0, \mu > 0$ due to the map's odd symmetry.

A sufficient condition for a trajectory beginning at $x_1 = \mu$ to be bounded is that the upper envelope $x_{n+1} \leq \mu + x_n^\rho$ intersect the identity line; that is, $\beta^\rho - \alpha + \mu = 0$ has a solution $\beta > 0$.
Noting that $F(x) = x^\rho - x + \mu$ has a minimum of $\rho^{\frac{1}{1-\rho}}\left(\frac1\rho-1\right)$ and that $F(0) = \mu > 0$, one sees that such a solution $\beta$ exists if $\mu \leq \rho^{\frac{1}{1-\rho}}\left(1-\frac1\rho\right)$; this region of parameter space corresponds to the region bounded by the green curve $\gamma_g$ in fig.~\ref{fig:chaosmeasures}A.
Evidenced by the existence of a positive Lyapunov exponent within this region, these bounded trajectories can nevertheless behave chaotically.
We now seek to prove that a trajectory $x_n$ with $x_1 = \mu$ converges to a stable fixed point when the map is an expansion ($\rho > 1$) and the splitting parameter $\mu$ is small.

One method to guarantee that a trajectory beginning at $x_1 = \mu$ converges to a fixed point is to establish a bound $x_n \leq \beta$ as before, subject to the additional constraint that $\left|\frac{\text{d}x_{n+1}}{\text{d}x_n}(x_n)\right| < 1$ for all $0 < x_n < \beta$.
Using the Brouwer fixed point theorem alongside the established bounds on the map's derivative, the existence of a unique fixed point of the map $x^*$ is verified within the interval $0 < x^* \leq \beta$ as $x_{n+1}(x_n) - x_n$ is monotone decreasing for $0 < x_n \leq \beta$.
Furthermore, this fixed point is determined to be stable.
In order to determine a large value $\mu$ such that a suitable $\beta$ exists, note that $\left|\frac{\text{d}x_{n+1}}{\text{d}x_n}\right| \leq x_n^{\rho-1}\sqrt{\rho^2+\omega^2}$: it is enough to satisfy $x_n^{\rho-1}\sqrt{\rho^2+\omega^2} < 1$ by choosing $\beta$ such that $x_n \leq \beta < (\rho^2+\omega^2)^{\frac{1}{2(1-\rho)}}$.
As $F(x)$ has its smallest positive root at $x = \beta$ and is a convex function, we can obtain an upper bound on $\beta$ by Jensen's inequality applied via a chord through $(0, F(0)) = (0, \mu)$ and $f$'s minimum $\left(\rho^{\frac{1}{1-\rho}}, \rho^{\frac{1}{1-\rho}}\left(\frac1\rho-1\right)+\mu\right)$: certainly $\beta \leq \frac{\mu}{1-\frac1\rho}$.
Hence a suitable bound $x_n \leq \beta$ exists if $\mu < \left(1-\frac1\rho\right)(\rho^2+\omega^2)^{\frac{1}{2(1-\rho)}}$; equality here yields the purple curve $\gamma_p$ in fig.~\ref{fig:chaosmeasures}A.
It is easy to see by the symmetry of the map that these stability conditions are nearly identical if $\mu < 0$; one needs only consider establishing the same bounds instead on the absolute value of $\mu$.

In the case of $\rho < 1$, the one-sided envelopes are convex and thus an invariant interval containing $x=\mu$ always exists.
An upper bound $\beta$ on trajectories in this case is given by the sufficient constraint $|x_n| \leq (|\mu|+1)^{\frac{1}{1-\rho}} \leq \beta$.

\subsection{Shrimp tails and symbolic robustness}

\begin{figure}[t]
\centering
\includegraphics[width=1.0\linewidth]{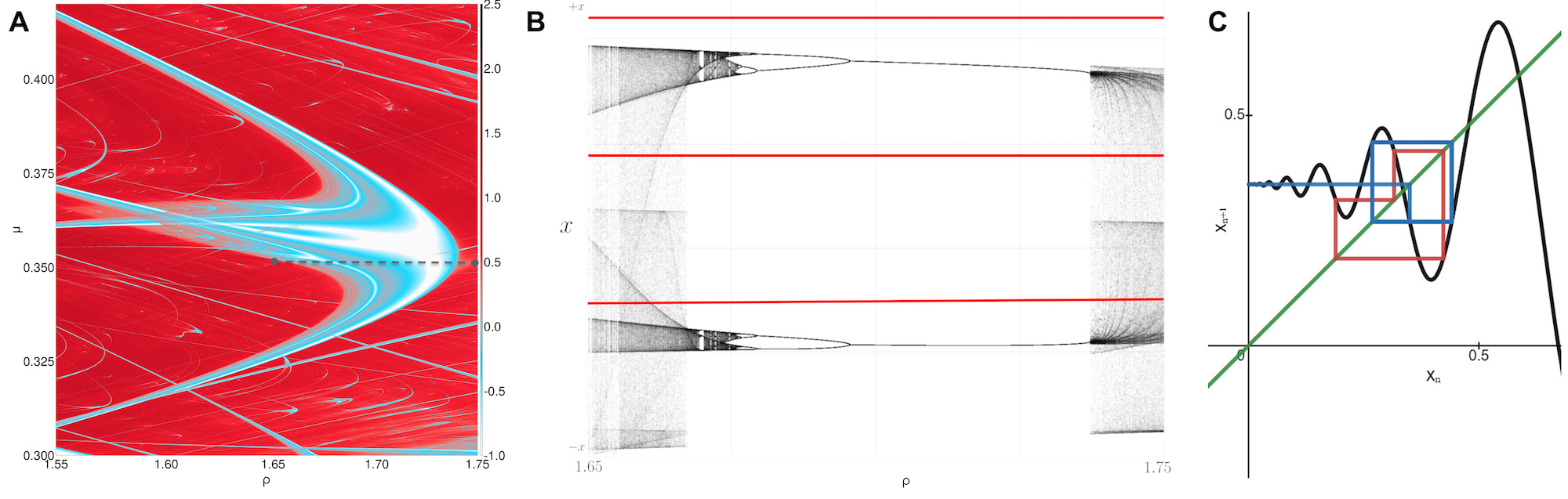}
\caption{(A) A two-dimensional Lyapunov-exponent sweep highlighting a ``shrimp'' structure, indicative of a saddle-node bifurcation in the one-dimensional saddle-focus map with $\omega = 10$. This region contains cascades of period doubling bifurcations of stable orbits of minimal period progressing in agreement with the Sharkovsky ordering\cite{shark3}. (B) The orbit diagram of a horizontal slice through the shrimp from (A) at $\mu=0.35$ (dotted interval). Observe the occurrence of a period doubling cascade in decreasing $\rho$ and a subsequent progression through orbits of periods with odd factors, signifying alongside the negative Lyapunov exponents in panel A that trajectories near the saddle-focus exhibit stable behavior. Nevertheless, there still exist robust orbits of period-3 in intervals containing $\mu$ throughout the shrimp. Such a period-3 orbit has been overlaid in red. Varying $\rho$ throughout the shrimp continuously deforms this red orbit, preserving its 3-periodicity.
(C) Juxtaposition of cobwebs of a period-3 orbit and a $x_1 = \mu$ orbit (from the saddle-focus) within the shrimp from panel A at a choice of $\rho$ where $\mu$ converges to a period-2 orbit.}\label{fig:shrimporbitdiagram}
\end{figure} 

Figure \ref{fig:shrimporbitdiagram} presents a detailed exploration of a ``shrimp'' structure identified in the Lyapunov-exponent scan (fig.~\ref{fig:shrimporbitdiagram}A) of the one-dimensional map with $\omega = 10$.
These regions arise from saddle-node bifurcations and exhibit periodic orbits robust to perturbations both in parameter space and in the one-dimensional interval map (\ref{eqn:1dmaps}).
One key observation is the presence of period doubling cascades, a common indicator of the emergence of chaotic dynamics.
Moreover, from the orbit diagram in figure \ref{fig:shrimporbitdiagram}B we observe that the periods of these orbits appear to progress monotonically through Sharkovsky's order\cite{shark3}.
The ``tails'' of these shrimp, those long negative-Lyapunov-exponent regions along decreasing $\rho$, carry on all the way to $\rho = 0$ and beyond, though the shrimp may be partially obscured by multistability.
The existence of these features, keeping multistability in mind, expands our understanding of the chaotic nature of the saddle-focus map and sets the stage for more in-depth study.

Figure~\ref{fig:shrimporbitdiagram}B depicts on the vertical axis the branches of stable periodic orbits originating at the boundary of the shrimp, plotted against the bifurcation parameter $\rho$ at fixed $\mu = 0.35$.
These periodic orbits develop in a manner reminiscent of saddle-node bifurcations and their further development in unimodal maps.
Despite the presence of stable periodic orbits within the shrimp appearing to coincide in evolution of periodicities with the Sharkovsky order as $\rho$ decreases, period-3 orbits can be easily identified throughout the windows, as is seen in the juxtaposed red curves in fig.~\ref{fig:shrimporbitdiagram}B corresponding to a persistent period-3 orbit; a cobweb diagram of another period-3 orbit within the shrimp is depicted in fig.~\ref{fig:shrimporbitdiagram}C.
This is important to keep in mind going forward, as the existence of a period-doubling cascade and subsequent progression to odd-period cycles does not by the Sharkovsky theorem imply the nonexistence of period-3 orbits.
At the same time, the existence of the negative-Lyapunov-exponent shrimp structure tells one nothing about the existence or absence of chaotic sets within intervals bounded by period-two orbits; multistability is prevalent throughout saddle-focus systems.

Our computations of the Lempel-Ziv complexity\cite{lempel1976complexity} for a symbolic sequence at each parameter value in the $(\rho, \mu)$-plane are showcased in fig.~\ref{fig:chaosmeasures}B.
The Lempel-Ziv complexity is a measure of the complexity of binary sequences, related in purpose to the notion of Kolmogorov complexity; it is defined as the length of a partition of a finite binary sequence such that each element of the partition is the shortest substring not having already occurred, less the final element if it happens to be a duplicate.
For instance, the binary sequence [$010110010111$] is partitioned as $\{0, 1, 01, 10, 010, 11, 1\}$, so it has a Lempel-Ziv complexity of $6$.
After computing the Lempel-Ziv complexity $C$ of a symbolic sequence of length $N$, we normalize by taking $\overline{C} = \frac{\ln(N)}{N} C$, as is done in our recent publication~\cite{scully2021measuring}.
The region confined by the purple curve in the two-dimensional LZ-sweep, as shown in fig.\ref{fig:chaosmeasures}B, displays sequences of minimal complexity, with quick convergence to unique fixed points for $\mu \geq 0$ or period-2 orbits for $\mu < 0$. 
However, substantial regions associated with positive Lyapunov exponents (refer to fig.\ref{fig:chaosmeasures}A) similarly exhibits low symbolic complexity. Chaos here does not change sign, and thus does not interact with homoclinics.

Within the region populated by homoclinics in the complexity scan, there is a ``sheet'' of high complexity interspersed by stability windows.
These windows align with tails of shrimp structures visible in the Lyapunov-exponent scan in fig.\ref{fig:chaosmeasures}A. 
The sheet appears as noise, seemingly induced by the sensitivity of symbolic sequences to perturbations of their generating trajectories, while the windows of reduced complexity indicate robust convergence to specific symbolic sequences. 
Furthermore, the geometric organization of the level sets of very small symbolic complexities within these shrimp tails -- and also across much of the boundary of the region of nontrivial symbolic complexity -- mirrors that of the homoclinic curves seen in the symbolic sequence scans from fig.\ref{fig:fullhomoclinicstructure} due to transients.

\section{Conclusions and future directions}

In this study, we delved into the heart of chaos, exploring the rich dynamics inherent in low-dimensional systems of ODEs, particularly the map associated with the Shilnikov saddle-focus homoclinic bifurcation. 
Inspired by the foundational work of Sharkovsky in one-dimensional maps, our research adopted two primary approaches: the generation of binary sequences to symbolically represent the dynamical behavior at each point in the parameter space, and the subsequent geometric analysis of the one-dimensional map (\ref{eqn:1dmaps}) in elucidating homoclinic bifurcation structure.
These techniques unveiled the intricate details of the homoclinic bifurcation structures relating to the saddle-focus, shedding light on the complex organization of these orbits.

Utilizing the Lyapunov exponent enabled us to illustrate the chaotic regions and stability zones within the saddle-focus map's parameter plane.
However, this method falls short when addressing multistability.
Our research revealed that the stability region of the saddle-focus map dramatically narrows near the codimension-two point $(\mu = 0$, $\rho = 1)$, representative of the Belyakov case~\cite{Belyakov3}.
This discovery raises profound questions about the nature of chaos at nonzero $\mu$, particularly as the $\rho = 1$ case corresponds to a nonhyperbolic saddle-focus, delicately balanced between the map's expansive and contractive behaviors.
Moreover, the relationship between the one-dimensional saddle-focus map and the corresponding two-dimensional return map has nuances that may result in obscuring chaotic behavior in the full saddle-focus ODE system.
The exploration of this theoretical frontier warrants deeper examination, and the one-dimensional map framework presents a promising avenue for this future endeavor, further building upon the pioneering work of L. P.~Shilnikov in the study of two- and higher-dimensional return maps.

Beyond this, there are additional aspects of both the stability and homoclinic structure that await scrutiny.
The occurrence of multistability within the map and its relationship to periodic orbits, as well as their corresponding homoclinics in systems of ODEs featuring a saddle-focus, represent fertile ground for future investigation.
In future research on these topics, we would like to:
\begin{itemize}
\item produce tools to extend our Lyapunov-exponent scans along stability branches,
\item visualize 2-dimensional homoclinic submanifolds of the $(\rho, \mu, \omega)$-parameter space by a method similar to the symbolic method we showcase in this paper, followed by an investigation of the homotopy types of these submanifolds, and
\item develop a computational method for efficiently scanning the $(\rho, \mu)$-parameter plane for the Sharkovsky-largest minimal-period orbit exhibited at each parameter choice, demonstrating the level of periodicity within the Sharkovsky order.
\end{itemize}
Investigation into these areas will not only enhance our understanding of the rich dynamics in such systems but also contribute to the broader theoretical framework for analyzing complex dynamical systems.

In conclusion, our research stands as a testament to the enduring impact of Sharkovsky's groundbreaking work \cite{shark1,shark2,shark4,shark3} on one-dimensional maps.
Our methods, influenced by his research, not only simplify the analysis of intricate dynamical structures but also offer a promising avenue for future investigations into similar low-dimensional systems.
The broad applicability of these techniques makes a significant contribution to the mathematical toolbox for studying complex dynamics, underscoring their potential to advance our understanding of chaos and complex dynamical systems.

\section*{Acknowledgments}
We thank the Brains \& Behavior initiative of Georgia State University for the B\&B graduate fellowship awarded to J.~Scully. 

\section*{References}
\bibliography{sharkovsky_submission.bib}

\end{document}